\documentclass[12pt]{amsart}

\usepackage{graphicx}
\usepackage{amsmath}
\usepackage{hyperref}

\usepackage{tikz}

\def\frak{\mathfrak}

\def\R{\mathbb{R}}

\def\cD{\mathcal{D}}

\def\al{\alpha}

\def\si{\sigma}

\def\om{\omega}

\def\Om{\Omega}

\newcommand{\der}{{\rm d}}

\numberwithin{equation}{section}

\newtheorem{theorem}{Theorem}[section]

\newtheorem{proposition}[theorem]{Proposition}

\theoremstyle{remark}

\theoremstyle{remark}

\usepackage{amssymb,stmaryrd}
\usepackage{amscd}

\newcommand{\qr}{
\begin{center}
\includegraphics[scale=0.5]{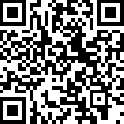}\\

\text{Scan the QR code to view more articles from the author}
\end{center}}

\makeatletter
\@namedef{subjclassname@2020}{%
  \textup{2020} Mathematics Subject Classification}
\makeatother

\author{Matthew Randall}
\address{Institute of Mathematical Sciences \\
ShanghaiTech University\\
393 Middle Huaxia Road\\
Shanghai 201210\\
China}
\email{mjrandall@shanghaitech.edu.cn}

\title{Local equivalence of some maximally symmetric $(2,3,5)$-distributions}


\begin{document}

\begin{abstract}
Using a complex parametrisation of $su(2)$, we show a change of coordinates that maps the maximally symmetric rolling $(2,3,5)$-distribution to the flat Cartan distribution. This establishes the local equivalence between the maximally symmetric rolling model and the flat Cartan or Hilbert-Cartan distribution. For the maximally symmetric rolling distribution, we write down the vector fields that bracket-generate to give the split real form of the Lie algebra of $\frak{g}_2$, with two of the vector fields in the bracket-generating set given by the span of the rolling distribution.
\end{abstract}

\maketitle

\pagestyle{myheadings}
\markboth{Randall}{Local equivalence of some maximally symmetric $(2,3,5)$-distributions}

\section{Introduction}

Let $\cD$ be a maximally non-integrable rank 2 distribution on a 5-manifold $M$. The maximally non-integrable condition of $\cD$ determines a filtration of the tangent bundle $TM$ given by
\[
\cD \subset [\cD,\cD] \subset [\cD,[\cD,\cD]]\cong TM.
\]
The distribution $[\cD, \cD]$ has rank 3 while the full tangent space $TM$ has rank 5, hence such a geometry is also known as a $(2,3,5)$-distribution. Let $M_{xyzpq}$ denote the 5-dimensional mixed order jet space $J^{2,0}(\R,\R^2) \cong J^2(\R,\R)\times \R$ with local coordinates given by $(x,y,z,p,q)=(x,y,z,y',y'')$ (see also \cite{tw13}, \cite{tw14}). Let $\cD_{F(x,y,z,y',y'')}$ denote the maximally non-integrable rank 2 distribution on $M_{xyzpq}$ associated to the underdetermined differential equation $z'=F(x,y,z,y',y'')$. This means that the distribution is annihilated by the following three 1-forms
\begin{align*}
\om_1=\der y-p \der x, \qquad \om_2=\der p-q \der x, \qquad  \om_3=\der z-F(x,y,z,p,q) \der x.
\end{align*}
Such a distribution $\cD_{F(x,y,z,y',y'')}$ is said to be in Monge normal form (see page 90 of \cite{tw13}). 
When $F(x,y,z,p,q)=q^2$, we obtain the 1-forms of the Hilbert-Cartan distribution. In Section 5 of \cite{conf}, it is shown how to associate canonically to such a $(2,3,5)$-distribution a conformal class of metrics of split signature $(2,3)$ (henceforth known as Nurowski's conformal structure or Nurowski's conformal metrics) such that the rank 2 distribution is isotropic with respect to any metric in the conformal class. The method of equivalence \cite{cartan1910} (also see the introduction to \cite{annur}, Section 5 of \cite{conf} and \cite{Strazzullo}) produces the 1-forms $(\theta_1, \theta_2,\theta_3, \theta_4, \theta_5)$ that gives a coframing for Nurowski's metric. These 1-forms satisfy the structure equations
\begin{align}\label{cse}
\der \theta_1&=\theta_1\wedge (2\Om_1+\Om_4)+\theta_2\wedge \Om_2+\theta_3 \wedge \theta_4,\nonumber\\
\der \theta_2&=\theta_1\wedge\Om_3+\theta_2\wedge (\Om_1+2\Om_4)+\theta_3 \wedge \theta_5,\nonumber\\
\der \theta_3&=\theta_1\wedge\Om_5+\theta_2\wedge\Om_6+\theta_3\wedge (\Om_1+\Om_4)+\theta_4 \wedge \theta_5,\\
\der \theta_4&=\theta_1\wedge\Om_7+\frac{4}{3}\theta_3\wedge\Om_6+\theta_4\wedge \Om_1+\theta_5 \wedge \Om_2,\nonumber\\
\der \theta_5&=\theta_2\wedge \Om_7-\frac{4}{3}\theta_3\wedge \Om_5+\theta_4\wedge\Om_3+\theta_5\wedge \Om_4,\nonumber
\end{align}
where $(\Om_1, \ldots, \Om_7)$ and two additional 1-forms $(\Om_8, \Om_9)$ together define a rank 14 principal bundle over the 5-manifold $M$ (see \cite{cartan1910} and Section 5 of \cite{conf}). A representative metric in Nurowski's conformal class \cite{conf} is given by
\begin{align}\label{metric}
g=2 \theta_1 \theta_5-2\theta_2 \theta_4+\frac{4}{3}\theta_3 \theta_3.
\end{align}
When $g$ has vanishing Weyl tensor, the distribution is called maximally symmetric and has split $G_2$ as its group of local symmetries. For further details, see the introduction to \cite{annur} and Section 5 of \cite{conf}. For further discussion on the relationship between maximally symmetric $(2,3,5)$-distributions and the automorphism group of the split octonions, see Section 2 of \cite{tw13}.

$(2,3,5)$-distributions also arise from the study of the configuration space of two surfaces rolling without slipping or twisting over each other \cite{AN14}, \cite{BH} and \cite{BM}. The configuration space can be realised as the An-Nurowski circle twistor distribution \cite{AN14} and in the case of two spheres with radii in the ratio $1:3$ rolling without slipping or twisting over each other, there is again maximal $G_2$ symmetry. 

In the work of \cite{r17}, a description of maximally symmetric $(2,3,5)$-distributions obtained from Pfaffian systems with $SU(2)$ symmetry was discussed and its relationship with the rolling distribution was investigated. In particular, the An-Nurowski circle twistor bundle can be realised by considering the Riemannian surface element of the unit sphere arising from one copy of $SU(2)$ and the other Riemannian surface element with Gaussian curvature $9$ or $\frac{1}{9}$ from another copy of $SU(2)$. Both Lie algebras of $su(2)$ are parametrised by the left-invariant vector fields. See \cite{r17} for further details. 

In this article we accomplish three things: 1) Instead of using the left-invariant vector fields, by using a complex parametrisation of $su(2)$ analogous to that done in \cite{r21a} whilst retaining the {\it rolling} terminology for this complex case, we show a change of coordinates that brings the 1-forms annihilating the rolling distribution to the 1-forms encoding the Monge equation
\begin{align}\label{meq}
\der y-p \der x,\nonumber\quad \der p-q \der x,\quad \der z-\left(qz^2-\frac{1}{\al^2+1}(\sqrt{q}z-\frac{1}{2\sqrt{q}x})^2\right) \der x,
\end{align} 
which was originally derived in \cite{r21a}. In other words we show that this complex version of the sphere rolling distribution can be brought into the Monge normal form with 
\[
F(x,y,z,p,q)=\left(qz^2-\frac{1}{\al^2+1}(\sqrt{q}z-\frac{1}{2\sqrt{q}x})^2\right). 
\]
Here $\al$ is a complex number, and the maximally symmetric case is obtained 
whenever $\al^2=-\frac{1}{9}$ or $\al^2=-9$.

2) In the aforementioned maximally symmetric case, we find the change of coordinates that maps the rolling distribution into the flat Cartan distribution, and therefore as a corollary into the Hilbert-Cartan distribution. This establishes the local equivalence between the maximally symmetric rolling model and the flat Cartan or Hilbert-Cartan distribution. 

3) For the maximally symmetric rolling distribution, we write down the vector fields that bracket-generate to give the split real form of the Lie algebra of $\frak{g}_2$, with two of the vector fields in the bracket-generating set given by the span of the rolling distribution. These are presented in Theorems \ref{g2a} and \ref{g2b}

The computations here are done utilising heavily the \texttt{DifferentialGeometry} package in MAPLE 2018. 

\section{Flat Cartan distribution}
The coframe data of the canonical maximally symmetric $(2,3,5)$-distribution is given by 
\begin{align*}
\der\theta_1=\theta_3 \wedge \theta_4, \quad
\der\theta_2=\theta_3 \wedge \theta_5,\quad
\der\theta_3=\theta_4 \wedge \theta_5,\quad
\der\theta_4=0,\quad
\der\theta_5=0.
\end{align*}
This is the historic case studied by Cartan (\cite{cartan1893}, \cite{cartan1910}) and Engel (\cite{Engel}, \cite{Engel2}). 
There are local coordinates $(a_1,a_2,a_3,a_4,a_5)$ (see pages 159--160 of \cite{cartan1910}) such that 
\begin{align}\label{t1}
\theta_1&=\der a_1+\left(a_3+\frac{1}{2}a_4 a_5\right)\der a_4,\\
\label{t2}
\theta_2&=\der a_2+\left(a_3-\frac{1}{2}a_4 a_5\right)\der a_5,\\
\label{t3}
\theta_3&=\der a_3+\frac{1}{2}a_4 \der a_5-\frac{1}{2}a_5 \der a_4,\\
\theta_4&=\der a_4,\nonumber\\
\theta_5&=\der a_5.\nonumber
\end{align}
For non-zero constant $k$, the weighted rescaling
\[
(\al_1,\al_2,\al_3, \al_4, \al_5)\mapsto (k^3\al_1,k^3\al_2,k^2 \al_3, k\al_4, k\al_5)
\]
preserves $\theta_1$, $\theta_2$ and $\theta_3$.
We shall refer to this distribution annhilated by the 1-forms $\{\theta_1, \theta_2, \theta_3\}$ in (\ref{t1}), (\ref{t2}), (\ref{t3}) as the flat Cartan distribution. To facilitate our writing of the vector fields that bracket-generate the Lie algebra of split ${\frak g}_2$, let us
pass to the 1-forms 
\begin{align*}
\Theta_1&=dc_1-2c_4dc_3-4c_3dc_4,\\
\Theta_2&=dc_2+2c_5dc_3+4c_3 dc_5,\\
\Theta_3&=dc_3+c_5dc_4-c_4dc_5,
\end{align*}
by taking the change of coordinates
\begin{align*}
(c_1,c_2,c_3,c_4,c_5)&=\left(6a_1-2a_3a_4+a_4^2a_5,6a_2-2a_3a_5-a_4a_5^2,2a_3,-a_4,a_5\right).
\end{align*}
It follows that
\begin{align*}
\Theta_1=6\theta_1+2a_4 \theta_3,\quad \Theta_2=6\theta_2+2a_5\theta_3,\quad \Theta_3=2\theta_3,
\end{align*}
so the 1-forms $\{\Theta_1, \Theta_2, \Theta_3\}$ are in the span of $\{\theta_1,\theta_2,\theta_3\}$.
If we take $\frak{r_1}=c_5$,$\frak{r_2}=c_4$, $\frak{r_3}=c_3$, $\frak{r_4}=-\frac{1}{2}(c_2+3c_3c_5)$, $\frak{r_5}=\frac{1}{2}(c_1-3c_3c_4)$, then the 1-forms 
\begin{align*}
\der\frak{r_3}+\frak{r_1}\der \frak{r_2}-\frak{r_2}\der \frak{r_1}&=\Theta_3,\\
\der\frak{r_4}+\frac{1}{2}(\frak{r_3}\der \frak{r_1}-\frak{r_1}\der \frak{r_3})&=-\frac{1}{2}\Theta_2,\\
\der\frak{r_5}+\frac{1}{2}(\frak{r_2}\der \frak{r_3}-\frak{r_3}\der \frak{r_2})&=\frac{1}{2}\Theta_1,
\end{align*}
obtained are those annihilating the flat Engel distribution as given in \cite{Engel} and \cite{Engel2}.

Let us write down the vector fields
\begin{align*}
Z^1=\partial_{c_3}+2c_5\partial_{c_2}-2c_4\partial_{c_1},\\
Z^2=\partial_{c_4}+4c_3\partial_{c_1}-2c_5\partial_{c_3},\\
Z^3=\partial_{c_5}+2c_4\partial_{c_3}-4c_3\partial_{c_2},
\end{align*}
and define
\begin{align}\label{cdist}
S^1=Z^2+c_5Z^1,\hspace{12pt} S^2=Z^3-c_4Z^1, \hspace{12pt} S^3=-c_1Z^2+c_2Z^3-(c_1c_5+c_2c_4+c_3^2)Z^1.
\end{align}
The vector fields $S^1$ and $S^2$ are in the span of the distribution and are annihilated by the 1-forms $\{\Theta_1, \Theta_2, \Theta_3\}$. 
We say that the vector fields $\{S^1,S^2,S^3\}$ pairwise bracket-generate the Lie algebra of split ${\frak g}_2$ if the following holds:
defining
\begin{align*}
S^4&=[S^1,S^2],\quad S^5=[S^2,S^3],\quad S^6=[S^3,S^1],\\
L^1&=[S^1,S^4],\quad L^3=[S^2,S^5],\quad L^5=[S^3,S^6],\\
L^2&=[S^2,S^4],\quad L^4=[S^3,S^5],\quad L^6=[S^1,S^6]
\end{align*}
and
\begin{align*}
H=[S^2,S^6], \quad h=[S^4,S^3],
\end{align*}
we require that the set of vector fields
\begin{align*}
\{S^1,S^2,S^3,S^4,S^5,S^6,\frac{1}{4}(h-H),\frac{\sqrt{3}}{12}(h+H),L^1,L^2,L^3,L^4,L^5,L^6\}
\end{align*}
form the 14-dimensional Lie algebra of split ${\frak g}_2$ with the Cartan subalgebra spanned by $\frac{1}{4}(h-H)$ and $\frac{\sqrt{3}}{12}(h+H)$ and the root diagram given by the picture below with respect to this choice of the Cartan subalgebra. 
\begin{figure}[h!]
\begin{tikzpicture}
	\draw [stealth-stealth](-1,0) -- (1,0);
\draw (1,0) node[anchor=west] {{\tiny $S^1$}};
\draw (-1,0) node[anchor=east] {{\tiny $S^5$}};
	\draw [stealth-stealth](0,-1.732) -- (0,1.732);
\draw (0,-1.732) node[anchor=north] {{\tiny $L^5$}};
\draw (0,1.732) node[anchor=south] {{\tiny $L^2$}};
\draw [stealth-stealth](-0.5,-0.866) -- (0.5,0.866);
\draw (-0.5,-0.866) node[anchor=north] {{\tiny $S^{3}$}};
\draw (0.5,0.866) node[anchor=south] {{\tiny $S^4$}};
\draw [stealth-stealth](-1.5,-0.866) -- (1.5,0.866);
\draw (-1.5,-0.866) node[anchor=north] {{\tiny $L^4$}};
\draw (1.5,0.866) node[anchor=south] {{\tiny $L^1$}};
\draw [stealth-stealth](1.5,-0.866) -- (-1.5,0.866);
\draw(1.5,-0.866) node[anchor=north] {{\tiny $L^6$}};
\draw (-1.5,0.866) node[anchor=south] {{\tiny $L^3$}};
\draw [stealth-stealth](0.5,-0.866) -- (-0.5,0.866);
\draw (0.5,-0.866) node[anchor=north] {{\tiny $S^{6}$}};
\draw (-0.5,0.866) node[anchor=south] {{\tiny $S^2$}};
\end{tikzpicture}
\end{figure}

\begin{proposition}
The vector fields given in (\ref{cdist}) pairwise bracket-generate the Lie algebra of split $\frak{g}_2$.
\end{proposition}
Since $S^1$ and $S^2$ are spanned by the distribution, which is already given as part of the data, the non-trivial part in determining the generating set for the Lie algebra of split $\frak{g}_2$ for a maximally symmetric $(2,3,5)$-distribution is to find $S^3$. The way to find it is outlined as above. We find the change of coordinates that brings the 1-forms annihilating the distribution to the span of $\theta_1$, $\theta_2$, $\theta_3$ in the flat Cartan distribution. Then we determine the functions $(c_1,c_2,c_3,c_4,c_5)$ and write down the vector fields $Z^1$, $Z^2$, $Z^3$, which now determine the bracket-generating set of vector fields $S^1$, $S^2$, $S^3$ completely. 
The Lie algebra of split $\frak{g}_2$ that arises in this way can be viewed also as the infinitesimal symmetry algebra of the $(2,3,5)$-distribution annihilated by the 1-forms $\{\tilde \theta_1, \tilde \theta_2,\tilde \theta_3\}$ where
\begin{align*}
\tilde \theta_1&=\der \tilde a_1+\left(\tilde a_3+\frac{1}{2}\tilde a_4 \tilde a_5\right)\der \tilde a_4,\\
\tilde \theta_2&=\der \tilde a_2+\left(\tilde a_3-\frac{1}{2}\tilde a_4 \tilde a_5\right)\der \tilde a_5,\\
\tilde \theta_3&=\der \tilde a_3+\frac{1}{2}\tilde a_4 \der \tilde a_5-\frac{1}{2}\tilde a_5 \der \tilde a_4,
\end{align*}
under the transformation
\begin{align*}
(\tilde a_1,\tilde a_2,\tilde a_3,\tilde a_4,\tilde a_5)=(a_1+a_3 a_4,a_2+a_3 a_5,-a_3,a_4,a_5).
\end{align*}
The vector fields in the symmetry algebra are precisely the ones for which Lie derivative of $\tilde \theta_1$, $\tilde \theta_2$, $\tilde \theta_3$ with respect to these vector fields are in the span of $\{\tilde \theta_1, \tilde \theta_2,\tilde \theta_3\}$. See \cite{tw14} for more explanation. The \texttt{InfinitesimalSymmetriesOfEDS} command in the \texttt{ExteriorDifferentialSystems} subpackage of MAPLE solves these equations for a given Pfaffian system specified by the 1-forms $\tilde \theta_1$, $\tilde \theta_2$, $\tilde \theta_3$ to obtain the vector fields in the symmetry algebra, and a key part of the results that appear here is brought about running this command and attributable to the visionary creators behind it. 

To illustrate the procedure for writing down the Lie algebra of split $\frak{g}_2$, let us look at the example of the Hilbert-Cartan distribution. To map the Hilbert-Cartan distribution given by the annihilator of the 1-forms
\[
\der y-p \der x,\quad \der p-q \der x,\quad \der z-q^2\der x,
\]
into the flat Cartan distribution, we take
\begin{align*}
a_1=2z+2q^2x-4pq,\quad a_2=2y,\quad a_3=2p-qx,\quad a_4=2q,\quad a_5=-x.
\end{align*}
This gives
\begin{align*}
\theta_1=2(\der z-q^2 \der x)-4q(\der p-q \der x),\quad \theta_2=2(\der y-p \der x), \quad \theta_3=2(\der p-q \der x).
\end{align*}
We determine
\begin{align*}
c_1&=12z-32pq+12q^2 x,\quad c_2=12y+4p x-4q^2x,\\
c_3&=4p-2qx,\quad c_4=-2q,\quad c_5=-x.
\end{align*}
Finding the basis of vector fields $\partial_{c_1}$, $\partial_{c_2}$, $\partial_{c_3}$, $\partial_{c_4}$, $\partial_{c_5}$, we obtain
\begin{align*}
Z^1&=\frac{1}{4}(\partial_p-x\partial_y+4q\partial_z),\\
Z^2&=\frac{x}{4}\left(\partial_p-x\partial_y+4q\partial_z\right)-\frac{1}{2}\partial_q,\\
Z^3&=-(\partial_x+p\partial_y+q^2\partial_z+q\partial_p)-\frac{q}{2}(\partial_p-x\partial_y+4q\partial_z),
\end{align*}
so that
\begin{align*}
S^1&=-\frac{1}{2}\partial_q,\\
S^2&=-(\partial_x+p \partial_y+q^2\partial_z+q\partial_p),\\
S^3&=\frac{1}{2}(12z-32pq+12 q^2x)\partial_q-(12y+4 px-4 qx^2)(\partial_x+p \partial_y+q^2\partial_z+q\partial_p)\\
&\quad{}-(2p-qx)^2(\partial_p-x\partial_y+4q\partial_z).
\end{align*}
The vector fields $S^1$ and $S^2$ are in the span of the Hilbert-Cartan distribution and together with $S^3$ they pairwise bracket-generate to form a split $\frak{g}_2$ Lie algebra. 
In this paper we compute the bracket-generating set of vector fields for a complex version of the maximally symmetric spheres rolling distribution. In a sequel to this paper, we do the same calculation for the maximally symmetric distribution determined by solutions of the $k=\frac{2}{3}$ and $\frac{3}{2}$ generalised Chazy equation as discussed in \cite{r16} and \cite{r19}, and the calculation for the real hyperboloids rolling distribution in \cite{r21a}. Together they establish the equivalences of the maximally symmetric rolling distribution, maximally symmetric generalised Chazy distribution and the Hilbert-Cartan distribution to one another. 

\section{A rolling distribution and its Monge normal form}
We now consider the rolling distribution (with the coordinates slightly changed) given in Section 2 of \cite{r17}. The 1-forms that annihilate the $(2,3,5)$-distribution given by the span of 
\begin{align*}
X^1&=\partial_x-E^2,\\
X^2&=i\partial_q-e^{\al x}E^1+\al e^{\al x}E^3
\end{align*}
are given by 
\begin{align}\label{rolling1form}
\om_1=\si_1-i\exp(\al x)\der q, \quad \om_2=\si_2+\der x, \quad \om_3=\si_3+i\al \exp(\al x) \der q. 
\end{align}
It was shown in \cite{r17} using the left-invariant 1-forms on $SU(2)$ that the distribution encoded by these three 1-forms can be related to the distribution associated to the rolling of two spheres over each other without slipping or twisting, and that the distribution is maximally symmetric when $\al^2=-\frac{1}{9}$ or $\al^2=-9$. Depending on the sign of $\al=\pm \frac{i}{3}$ or $\al=\pm 3i$, we distinguish between the vector fields $X^{2+}$ and $X^{2-}$.
Let us now find the change of coordinates that maps it into the Monge normal form found in \cite{r21a}.

Instead of using the real left-invariant 1-forms on $SU(2)$, we use the complex parametrisation given by
\begin{align}\label{su2i}
\si_1&=i(\der y+y\der z)-i(y^2-1)(\der p-p \der z),\nonumber\\
\si_2&=\der y+y\der z-(y^2+1)(\der p-p \der z), \\
\si_3&=-2iy(\der p-p \der z)+i \der z, \nonumber
\end{align} 
where $i^2=-1$. However we still retain the terminology of refering to this distribution as rolling. If we feel pedantic, then we can make the distinction between real rolling when the lie algebra $su(2)$ is parametrised by the real left-invariant vector fields and {\it complex rolling} when $su(2)$ is parametrised by vector fields dual to the complex 1-forms (\ref{su2i}) above. The latter terminology is not apt since rolling should be a real physical process, but not entirely unjustified either following \cite{BH}. We make some remarks at the end of the paper about how to relate the real rolling distribution from this complex one. Nurowski's conformal structure for the complex distribution annihilated by the 1-forms (\ref{rolling1form}) given here can be computed without much change from \cite{r17}.

It can be verified that these three 1-forms satisfy the relations 
\[
\der \si_1=\si_2 \wedge \si_3, \quad \der \si_2=\si_3 \wedge \si_1, \quad \der \si_3=\si_1 \wedge \si_2.
\]
The vector fields dual to the 1-forms are given by
\begin{align*}
E^1&=\frac{i}{2}(y^2-1)\partial_y-\frac{i}{2}(2yp+1)\partial_p-i y\partial_z,\\
E^2&=\frac{1}{2}(y^2+1)\partial_y-\frac{1}{2}(2yp+1)\partial_p-y\partial_z,\\
E^3&=i(y\partial_y-p\partial_p-\partial_z).
\end{align*}
They satisfy the $su(2)$ Lie algebra
\[
[E^1,E^2]=-E^3,\quad [E^3,E^1]=-E^2,\quad [E^2,E^3]=-E^1.
\]
With this choice of vector fields we have the complex conjugates $\overline{E^1}=-E^1$, $\overline{E^3}=-E^3$ and so the complex conjugate of $X^{2+}$ gives $\overline{X^{2+}}=-X^{2-}$.
Let us define the 1-forms $t_1$, $t_2$ and $t_3$ given by the combinations
\begin{align*}
t_1&=\al \om_1+\om_3\\
&=i\al \der y+i((\al y^2+2y-\al)p+(\al y+1))\der z-i(\al y^2+2y-\al)\der p,\\
t_2&=\om_1-i \om_2=-i\der x-2ip \der z+2i \der p-i\exp(\al x) \der q,\\
t_3&=\om_3+y(\om_1-i \om_2)=-i y \der x+i \der z-i\exp(\al x) (y-\al) \der q.
\end{align*}

Making the change of coordinates
\begin{align*}
\hat y&=\exp(z)(\al y+1),\quad \hat p=-\exp(-z) p,\quad \hat x=q-\frac{1}{\al}\exp(-\al x),\\
\hat z&=\exp(z-\al x),\quad \hat q=\exp(-\al x), 
\end{align*}
we obtain
\begin{align*}
\frac{i}{2}\exp(-z)t_2&=\der \hat p+\frac{1}{2\hat z}\der \hat x, \\
-i\exp(z-\al x) t_3&=\der \hat z-\exp(z) (y-\al)\der \hat x,\\
-i \exp(z) t_1&=\der \hat y+\exp(2z)(\al y^2+2y-\al)\der \hat p.
\end{align*}
Now we have
\begin{align*}
\exp(z)=\frac{\hat z}{\hat q},\quad \al y+1=\frac{\hat q \hat y}{\hat z},\quad  y=\frac{1}{\al}\left(\frac{\hat q \hat y}{\hat z}-1\right).
\end{align*}
We therefore get
\begin{align*}
\frac{i}{2}\exp(-z)t_2&=\der \hat p+\frac{1}{2\hat z}\der \hat x, \\
-i\exp(z-\al x) t_3&=\der \hat z-\frac{\hat z}{\hat q}\left(\frac{1}{\al}\big(\frac{\hat q \hat y}{\hat z}-1\big)-\al\right)\der \hat x\\
&=\der \hat z-\left(\frac{1}{\al}\hat y-\big(\frac{1}{\al}+\al\big)\frac{\hat z}{\hat q}\right)\der \hat x,\\
-i \exp(z) t_1&=\der \hat y+\frac{\hat z^2}{\hat q^2}\left(\frac{1}{\al}\big(\frac{\hat q \hat y}{\hat z}-1\big)^2+\frac{2}{\al}\big(\frac{\hat q \hat y}{\hat z}-1\big)-\al\right)\der \hat p\\
&=\der \hat y+\frac{\hat z^2}{\hat q^2}\left(\frac{1}{\al}\frac{\hat q^2 \hat y^2}{\hat z^2}-\frac{1}{\al}-\al\right)\der \hat p\\
&=\der \hat y+\left(\frac{1}{\al}\hat y^2-\big(\frac{1}{\al}+\al\big)\frac{\hat z^2}{\hat q^2}\right)\der \hat p.
\end{align*}
Now take 
\[
\bar y=\hat x+2 \hat z \hat p,\quad \bar x=2\hat z,\quad \bar p=\hat p,\quad \bar z=-\frac{\hat y}{\al},\quad \bar q=-\frac{\al \hat q}{4\hat z (\hat y\hat q-(\al^2+1)\hat z)}.
\]
If we let
\[
F(\bar x, \bar y, \bar z, \bar p, \bar q)=\bar q \bar z^2-\frac{1}{\al^2+1}(\sqrt{\bar q}\bar z-\frac{1}{2\sqrt{\bar q}\bar x})^2=-\frac{\hat y^2\hat q^2-(\al^2+1)\hat z^2}{4\al\hat q \hat z(\hat y \hat q-(1+\al^2)\hat z)},
\]
then we get
\begin{align*}
\der \bar y-\bar p\der \bar x&=\der \hat x+2 \hat z \der \hat p,\\
\der \bar p-\bar q \der \bar x&=\der \hat p+\frac{2\al \hat q}{4\hat z(\hat y\hat q-(\al^2+1)\hat z)}\der \hat z=-\frac{1}{2\hat z}\der\hat x+\frac{\al \hat q}{2\hat z(\hat y\hat q-(\al^2+1)\hat z)}\der \hat z\\
&=-\frac{1}{2\hat z}\left(\der \hat x-\frac{\al \hat q}{\hat y\hat q-(\al^2+1)\hat z}\der \hat z\right),\\
\der \bar z-F(\bar x, \bar y, \bar z, \bar p, \bar q)\der \bar x&=-\frac{1}{\al}\der \hat y+\frac{2(\hat y^2\hat q^2-(\al^2+1)\hat z^2)}{4\al\hat q \hat z(\hat y \hat q-(1+\al^2)\hat z)}\der \hat z\\
&=-\frac{1}{\al}\der \hat y+\frac{2(\hat y^2\hat q^2-(\al^2+1)\hat z^2)}{\al\hat q}\left(-\frac{1}{2\al \hat q}\der \hat p\right)\\
&=-\frac{1}{\al}\der \hat y-\frac{\hat y^2\hat q^2-(\al^2+1)\hat z^2}{\al^2\hat q^2}\der \hat p\\
&=-\frac{1}{\al}\left(\der \hat y+\big(\frac{1}{\al}\hat y^2-\frac{\al^2+1}{\al}\frac{\hat z^2}{\hat q^2}\big)\der \hat p\right).
\end{align*}
The right hand side are multiples of $t_1, t_2, t_3$, and therefore in the span of the 1-forms $\om_1, \om_2, \om_3$.
We obtain the following theorem.
\begin{theorem}
Let $\om_1$, $\om_2$, $\om_3$ be the 1-forms given in (\ref{rolling1form}) with the complex parametrisation of $su(2)$ given in (\ref{su2i}). By the change of coordinates
\begin{equation}\label{rtom}
(\bar{x}, \bar{y}, \bar{z}, \bar{p}, \bar{q})=\left(2 e^{z-\al x},q-\frac{1}{\al}e^{-\al x}-2p e^{-\al x}, -\frac{\al y+1}{\al}e^z, -e^{-z} p, -\frac{1}{4e^{2z-\al x}(y-\al)}\right),
\end{equation}
the 1-forms $\om_1$, $\om_2$, $\om_3$ can be brought into the Monge normal form
\[
\der \bar y-\bar p \der \bar x,\quad \der \bar p-\bar q \der \bar x,\quad \der \bar z-\left(\bar q \bar z^2-\frac{1}{\al^2+1}(\sqrt{\bar q}\bar z-\frac{1}{2\sqrt{\bar q}\bar x})^2\right)\der \bar x.
\]
The distribution is maximally symmetric whenever $\al^2=-\frac{1}{9}$ or $\al^2=-9$.
\end{theorem}

The fact that the distribution is maximally symmetric for the parameters given follows from the computations done in \cite{r17} and \cite{r21a}.
We now exhibit the map from the 1-forms encoding the maximally symmetric Monge equation above to the flat Cartan distribution, and as a consequent obtain a map from the maximally symmetric rolling distribution to the flat Cartan distribution. We also write down the corresponding vector fields that bracket-generate to give the split $\frak{g}_2$ Lie algebra. 

\section{Local equivalence of the maximally symmetric rolling distribution to flat Cartan distribution}

In the two maximally symmetric cases, we get
\[
F(x,y,z,p,q)=-\frac{1}{8}\bar q\bar z^2+\frac{9}{8}\frac{\bar z}{\bar x}-\frac{9}{32}\frac{1}{\bar x^2\bar q}
\]
when $\al^2=-\frac{1}{9}$ and 
\[
F(x,y,z,p,q)=\frac{9}{8}\bar q\bar z^2-\frac{1}{8}\frac{\bar z}{\bar x}+\frac{1}{32}\frac{1}{\bar x^2\bar q}
\]
when $\al^2=-9$.
Let us now analyse the case when $\al^2=-\frac{1}{9}$.
When
\[
F=-\frac{1}{8}\bar q\bar z^2+\frac{9}{8}\frac{\bar z}{\bar x}-\frac{9}{32}\frac{1}{\bar x^2\bar q},
\]
we can map the 1-forms
\[
\bar\om_1=\der \bar y-\bar p \der \bar x,\quad \bar \om_2= \der \bar p-\bar q \der \bar x,\quad \bar \om_3=\der \bar z+\left(\frac{1}{8}\bar q\bar z^2-\frac{9}{8}\frac{\bar z}{\bar x}+\frac{9}{32}\frac{1}{\bar x^2\bar q}\right)\der \bar x
\]
into the flat Cartan distribution as follows.

Define
\begin{align*}
a_1&=-\frac{\sqrt{\bar x}}{6(2\bar q\bar x\bar z-1)^3}\left(-9(2\bar q\bar x\bar z+3)(2\bar q\bar x\bar z-1)^2\bar y+16 \bar q \bar x\bar z((2\bar q\bar x\bar z+3)(2\bar q\bar x\bar z-3))\bar x\bar p-512 \bar q^3\bar x^4\bar z^2\right),\\
a_2&=-\frac{\sqrt{\bar x}}{6(2\bar q\bar x\bar z-1)^3}\bigg(2\bar x(2\bar q\bar x\bar z+3)(20\bar q^2\bar x^2\bar z^2-12\bar q\bar x\bar z+9)\bar p^2-32\bar q\bar x^2(2\bar q\bar x\bar z-3)(2\bar q\bar x\bar z+3)\bar p\\
&\quad{} -9(2\bar q \bar x\bar z-1)^2(2\bar p\bar q\bar x\bar z+3\bar p+8\bar q \bar x)\bar y-64\bar q^2\bar x^3(10\bar q\bar x\bar z-9)\bigg),\\
a_3&=\frac{3}{2(2\bar q\bar x\bar z-1)}((2\bar q\bar x\bar z-1)(\bar y-\bar p \bar x)+8\bar q\bar x^2),\\
a_4&=-\frac{\sqrt{\bar x}}{(2\bar q\bar x\bar z-1)}(2\bar q\bar x\bar z+3),\\
a_5&=-\frac{\sqrt{\bar x}}{(2\bar q\bar x\bar z-1)}((2\bar q \bar x \bar z+3)\bar p+8\bar q\bar x).
\end{align*}
We find for this set of functions, 
\begin{align*}
\theta_1&=\frac{3}{2}\frac{\sqrt{\bar x}}{(2\bar q\bar x\bar z-1)}(2\bar q \bar x \bar z+3)\bar \om_1-\frac{8}{3}\frac{\bar x^{5/2}\bar q\bar z(2\bar q\bar x\bar z-3)(2\bar q\bar x\bar z+3)}{(2\bar q\bar x\bar z-1)^3}\bar \om_2-\frac{256}{3}\frac{\bar x^{9/2}\bar q^3\bar z}{(2\bar q\bar x\bar z-1)^3}\bar \om_3,\\
\theta_2&=\frac{3}{2}\frac{\sqrt{\bar x}}{(2\bar q\bar x\bar z-1)}(2\bar p\bar q \bar x \bar z+3\bar p+8 \bar q \bar x)\bar \om_1-\frac{8}{3}\frac{\bar x^{5/2}\bar q(2\bar q\bar x\bar z-3)(2\bar q\bar x\bar z+3)(\bar p\bar z+1)}{(2\bar q\bar x\bar z-1)^3}\bar \om_2\\
&\quad{} -\frac{256}{3}\frac{\bar x^{9/2}\bar q^3(\bar p\bar z+1)}{(2\bar q\bar x\bar z-1)^3}\bar \om_3,\\
\theta_3&=\frac{3}{2}\bar \om_1-\frac{(4\bar q^2\bar x^2\bar z^2-12\bar q\bar x\bar z-3)\bar x}{(2\bar q\bar x\bar z-1)^2}\om_2-\frac{32{\bar q}^2{\bar x}^3}{(2\bar q\bar x\bar z-1)^2}\bar \om_3,
\end{align*}
so that $\theta_1$, $\theta_2$ and $\theta_3$ is in the span of $\bar \om_1$, $\bar \om_2$, $\bar \om_3$.

Composing with the map (\ref{rtom})
when $\al=\pm \frac{i}{3}$, we obtain the transformation of $\om_1$, $\om_2$, $\om_3$ in (\ref{rolling1form}) into the flat Cartan distribution.

To derive the vector fields that generate $\frak{g}_2$, we compute
\begin{align*}
(c_1,c_2,c_3,c_4,c_5)=\left(6a_1-2a_3a_4+a_4^2a_5,6a_2-2a_3a_5-a_4a_5^2,2a_3,-a_4,a_5\right)
\end{align*}
and find
\begin{align*}
c_1&=-\frac{48\sqrt{2}}{(\al^2+1)^2\al}e^{\frac{z}{2}}\bigg(e^{-\frac{3\al x}{2}}(\al(7\al^2+1)y-\frac{1}{4}(15\al^4-10\al^2-1))-q e^{-\frac{\al x}{2}}\al (\al^2+1)(\al y+\frac{1-3\al^2}{4})\bigg),\\
c_2&=\frac{48\sqrt{2}}{(\al^2+1)^2\al}e^{-\frac{z}{2}}\bigg(e^{-\frac{3\al x}{2}}(\al(7\al^2+1)(yp+1)-\frac{1}{4}(15\al^4-10\al^2-1)p)\\
&-qe^{-\frac{\al x}{2}}\al(\al^2+1)(\al +(\al y+\frac{1-3\al^2}{4})p)\bigg),\\
c_3&=-\frac{3(9\al^2+1)}{\al(\al^2+1)}e^{-\al x}+3q,\\
c_4&=\frac{4\sqrt{2}}{(\al^2+1)}e^{\frac{z}{2}-\frac{\al x}{2}}(\al y+\frac{1-3\al^2}{4}),\\
c_5&=\frac{4\sqrt{2}}{(\al^2+1)}e^{-\frac{z}{2}-\frac{\al x}{2}}(\al+(\al y+\frac{1-3\al^2}{4})p).
\end{align*}
The 1-forms given by $\Theta_1$, $\Theta_2$ and $\Theta_3$ are in the span of $\om_1$, $\om_2$ and $\om_3$ precisely when $\al^2+\frac{1}{9}=0$.
When $\al=\pm\frac{i}{3}$, we obtain 
\begin{align*}
c_1&=-\frac{27\sqrt{2}}{2}e^{\frac{z}{2}}(e^{\mp \frac{ix}{2}}(y\pm i)\mp \frac{4}{3}iqe^{\mp \frac{ix}{6}}(y\mp i)),\\
c_2&=\frac{27\sqrt{2}}{2}e^{-\frac{z}{2}}(e^{\mp \frac{ix}{2}}(1+(y\pm i)p)\mp \frac{4}{3}iqe^{\mp \frac{ix}{6}}(1+(y\mp i)p)),\\
c_3&=3q,\\
c_4&=\pm i \frac{3\sqrt{2}}{2}e^{\frac{z}{2}\mp\frac{ix}{6}}(y\mp i),\\
c_5&=\pm i \frac{3\sqrt{2}}{2}e^{-\frac{z}{2}\mp\frac{ix}{6}}(1+(y\mp i)p).
\end{align*}
We can further rescale $(c_1,c_2,c_3,c_4,c_5)\mapsto (\frac{2\sqrt{2}}{27}c_1,\frac{2\sqrt{2}}{27} c_2, \frac{2}{9}c_3, \frac{\sqrt{2}}{3}c_4, \frac{\sqrt{2}}{3}c_5)$ to obtain
\begin{align*}
c_1&=-2e^{\frac{z}{2}}(e^{\mp \frac{ix}{2}}(y\pm i)\mp \frac{4}{3}iqe^{\mp \frac{ix}{6}}(y\mp i)),\\
c_2&=2e^{-\frac{z}{2}}(e^{\mp \frac{ix}{2}}(1+(y\pm i)p)\mp \frac{4}{3}iqe^{\mp \frac{ix}{6}}(1+(y\mp i)p)),\\
c_3&=\frac{2}{3}q,\\
c_4&=\pm i e^{\frac{z}{2}\mp\frac{ix}{6}}(y\mp i),\\
c_5&=\pm ie^{-\frac{z}{2}\mp\frac{ix}{6}}(1+(y\mp i)p).
\end{align*}
This gives 
\begin{align*}
Z^1&=-\frac{3i}{2}X^{2\pm}\pm \frac{1}{2}e^{\pm\frac{ix}{3}}(5E^3\pm 3iE^1),\\
Z^2&=\pm\frac{3i}{4}e^{-\frac{z}{2}\pm\frac{ix}{6}}((y\pm i)p+1)X^1\mp3e^{-\frac{z}{2}\mp\frac{ix}{6}}((y\mp i)p+1)X^{2\pm}\\
&\quad{} -\frac{i}{2}e^{-\frac{z}{2}\pm\frac{ix}{6}}((y\mp i)p+1)(5 E^3\pm 3iE^1),\\
Z^3
&=\mp\frac{3i}{4}e^{\frac{z}{2}\pm\frac{ix}{6}}(y\pm i)X^1\pm3(y\mp i)e^{\frac{z}{2}\mp\frac{ix}{6}}X^{2\pm}+\frac{i}{2}(y\mp i)e^{\frac{z}{2}\pm \frac{ix}{6}}(5 E^3\pm 3iE^1)
\end{align*}
and we obtain
\begin{align}\label{s1a}
S^1&=\pm\frac{3i}{4}e^{-\frac{z}{2}\pm\frac{ix}{6}}((y\pm i)p+1)X^1\mp\frac{3}{2}e^{-\frac{z}{2}\mp\frac{ix}{6}}((y\mp i)p+1)X^{2\pm},\\\label{s2a}
S^2&=\mp\frac{3i}{4}e^{\frac{z}{2}\pm\frac{ix}{6}}(y\pm i)X^1\pm \frac{3}{2}e^{\frac{z}{2}\mp\frac{ix}{6}}(y\mp i)X^{2\pm}.
\end{align}
These two vector fields lie in the span of $X^1$ and $X^{2\pm}$ and together with 
\begin{align}\label{s3a}
S^3&=\mp4iqX^1+\frac{2}{3}i(q^2-9e^{\mp\frac{2ix}{3}})X^2\mp \frac{2}{9}q^2e^{\pm\frac{ix}{3}}(5E^3\pm 3iE^1)
\end{align} 
they bracket-generate the Lie algebra of split $\frak{g}_2$.

\begin{theorem}\label{g2a}
For the maximally symmetric distribution spanned by $X^1$, $X^{2\pm}$ when $\al=\pm \frac{i}{3}$ and annihilated by the 1-forms in (\ref{rolling1form}) with $(\si_1,\si_2,\si_3)$ given by (\ref{su2i}), the Lie algebra of split $\frak{g}_2$ is obtained from the pairwise bracket-generating set $\{S^1,S^2,S^3\}$ where $S^1$, $S^2$ and $S^3$ are given in (\ref{s1a}), (\ref{s2a}), (\ref{s3a}).
\end{theorem}

It can be checked that $S^{1-}$, $S^{2-}$, $S^{3-}$ are the complex conjugates of $S^{1+}$, $S^{2+}$, $S^{3+}$ respectively, depending on the sign of $\al$.

When 
\[
F=\frac{9}{8}qz^2-\frac{1}{8}\frac{z}{x}+\frac{1}{32qx^2},
\] which corresponds to the case where $\al^2=-9$, we can repeat the procedure above to derive the mapping into the flat Cartan distribution and write down the corresponding vector fields that bracket generate $\frak{g}_2$.
To map the distribution given by kernel of the 1-forms
\begin{align*}
\der y-p \der x, \quad \der p-q \der x, \quad \der z-\left(-\frac{1}{8}qz^2+\frac{9}{8}\frac{z}{x}-\frac{9}{32qx^2}\right)\der x
\end{align*}
into 
\begin{align*}
\der \tilde y-\tilde p \der \tilde x, \quad \der \tilde p-\tilde q \der \tilde x, \quad \der \tilde z-\left(\frac{9}{8}\tilde q\tilde z^2-\frac{1}{8}\frac{\tilde z}{\tilde x}+\frac{1}{32\tilde q\tilde x^2}\right)\der \tilde x, 
\end{align*}
we take either the map
\begin{align*}
\tilde x=\frac{1}{x},\quad \tilde y=\frac{y}{x}, \quad \tilde z=\frac{z}{9x},\quad \tilde p=y-px, \quad \tilde q=q x^3
\end{align*}
or 
\begin{align*}
\tilde x&=\frac{9(2qxz-1)^2}{(2qxz-9)^2x},\quad \tilde y=\frac{9(2qxz-1)^2y-8xp(2qxz-3)(2qxz+3)+256q^2x^3z}{(2qxz-9)^2x},\\
 \tilde z&=\frac{(2qxz-1)}{2x^2 q(2q xz-9)},\quad \tilde p=y-px+\frac{16 x^2q}{2qxz-1}, \quad \tilde q=\frac{(2qxz-9)^2x^2}{18(2qxz-1)^3z}.
\end{align*}
Both maps compose twice to give the identity. An interesting observation is that composing one map with the inverse of the other gives a symmetry of a distribution to itself which squares to the identity.
The map
\begin{align*}
\tilde{\tilde x}&=\frac{(2qxz-9)^2x}{9(2qxz-1)^2}, \quad \tilde{\tilde y}=y-\frac{8x(2qxz-3)(2qxz+3)}{9(2qxz-1)^2}p+\frac{256q^2x^3z}{9(2qxz-1)^2},\\
\tilde{\tilde z}&=\frac{2qxz-9}{2(2qxz-1)qx},\quad \tilde{\tilde p}=p-\frac{16xq}{2qxz-9},\quad \tilde{\tilde q}=\frac{81(2qxz-1)^3}{2zx(2qxz-9)^3},
\end{align*}
for instance take the 1-forms 
\begin{align*}
\der y-p \der x, \quad \der p-q \der x, \quad \der z-\left(-\frac{1}{8}qz^2+\frac{9}{8}\frac{z}{x}-\frac{9}{32qx^2}\right)\der x
\end{align*}
into themselves.

Composing the map from the $\al^2=-\frac{1}{9}$ Monge equation to the flat Cartan distribution
\[
(\bar x, \bar y, \bar z, \bar p, \bar q) \mapsto (a_1,a_2,a_3,a_4,a_5)
\]
given above with
\[
\bar x=\frac{1}{\tilde x},\quad \bar y=\frac{\tilde y}{\tilde x},\quad \bar z=9\frac{\tilde z}{\tilde x},\quad \bar p=\tilde y-\tilde p\tilde x,\quad \bar q=\tilde q\tilde x^3,
\]
followed by (\ref{rtom}) with $\al^2=-9$, the coordinate transformation from the rolling distribution with $\al^2=-9$ to the flat Cartan distribution can be found and the functions $(c_1,c_2,c_3,c_4,c_5)$ are given by
\begin{align*}
c_1&=-\frac{18\sqrt{2}}{(8\al y+\al^2+9)^2}e^{-\frac{3z}{2}+\frac{\al x}{2}}\left((4\al y-\al^2+3)(8\al y+\al^2+9)p+2\al (10\al y-\al^2+9)\right),\\
c_2&=\frac{18\sqrt{2}}{\al(8\al y+\al^2+9)^2}\bigg(e^{-\frac{3}{2}z-\frac{\al x}{2}}(\frac{1}{3}(12\al y+5\al^2+9)(8\al y+\al^2+9)p+2\al(10\al y+3\al^2+9))\\
&\quad{} -q e^{-\frac{3}{2}z-\frac{\al x}{2}}((4\al y-\al^2+3)(8\al y+\al^2+9)p+2\al (10 \al y-\al^2+9))\al \bigg),\\
c_3&=-\frac{3}{8\al y+\al^2+9}e^{-z}((8\al y+\al^2+9)p+4\al),\\
c_4&=\frac{3\sqrt{2}}{2}\frac{4\al y-\al^2+3}{8\al y+\al^2+9}e^{-\frac{z}{2}+\frac{\al x}{2}},\\
c_5&=\frac{3\sqrt{2}}{2\al (8 \al y+\al^2+9)}\left(e^{-\frac{z}{2}-\frac{\al x}{2}}(4\al y+\frac{5}{3}\al^2+3)+\al e^{-\frac{z}{2}-\frac{\al x}{2}} q(4\al y-\al^2+3)\right).
\end{align*}
Again the 1-forms given by $\Theta_1$, $\Theta_2$ and $\Theta_3$ are in the span of $\om_1$, $\om_2$ and $\om_3$ precisely when $\al^2+9=0$.
When $\al=\pm 3 i$, we have
\begin{align*}
c_1&=\frac{9\sqrt{2}}{2y^2}e^{-\frac{3}{2}z\pm \frac{3i x}{2}}\left(2y(y\mp i)p+\frac{1}{4}(5y\mp 3i)\right),\\
c_2&=-\frac{9\sqrt{2}}{2y^2}e^{-\frac{3z}{2}}\bigg( qe^{\pm\frac{3ix}{2}}(2y(y\mp i)p+\frac{1}{4}(5y\mp 3i))\pm  \frac{i}{3}e^{\mp \frac{3 i x}{2}} (2y (y\pm i)p+\frac{1}{4}(5 y\pm 3i))\bigg),\\
c_3&=-\frac{3}{2y}e^{-z}(2py+1),\\
c_4&=\frac{3\sqrt{2}}{4y}e^{-\frac{z}{2}\pm \frac{3i x}{2}}(y \mp i),\\
c_5&=-\frac{3\sqrt{2}}{4y}e^{-\frac{z}{2}}\left(qe^{\pm\frac{3i x}{2}}(y\mp i)\pm\frac{i}{3}e^{\mp\frac{3i x}{2}}(y\pm i)\right).
\end{align*}

We can further rescale $(c_1,c_2,c_3,c_4,c_5)\mapsto (-\frac{32}{27\sqrt{2}}c_1,-\frac{32}{27\sqrt{2}} c_2,\frac{8}{9}c_3, -\frac{4}{3\sqrt{2}}c_4, -\frac{4}{3\sqrt{2}}c_5)$ to obtain
\begin{align*}
c_1&=\frac{16}{3y^2}e^{-\frac{3}{2}z\pm \frac{3i x}{2}}\left(2y(y\mp i)p+\frac{1}{4}(5y\mp 3i)\right),\\
c_2&=\frac{16}{3y^2}e^{-\frac{3z}{2}}\bigg( qe^{\pm\frac{3ix}{2}}(2y(y\mp i)p+\frac{1}{4}(5y\mp 3i))\pm \frac{i}{3} e^{\mp \frac{3 i x}{2}}(2y (y\pm i)p+\frac{1}{4}(5 y\pm 3i))\bigg),\\
c_3&=-\frac{4}{3y}e^{-z}(2py+1),\\
c_4&=-\frac{1}{y}e^{-\frac{z}{2}\pm \frac{3i x}{2}}(y \mp i),\\
c_5&=\frac{1}{y}e^{-\frac{z}{2}}\left(qe^{\pm\frac{3i x}{2}}(y\mp i)\pm\frac{i}{3}e^{\mp\frac{3i x}{2}}(y\pm i)\right).
\end{align*}
Let us denote the vector field $X^3$ by
\[
X^3=i\frac{y^2-1}{2y}E^1+\frac{y^2+1}{2y}E^2+5 i E^3.
\] 
Computing the basis of vector fields $\partial_{c_1}$, $\partial_{c_2}$, $\partial_{c_3}$, $\partial_{c_4}$, $\partial_{c_5}$, we get
\begin{align*}
Z^1&=-\frac{3}{4}e^z(y^2+1) X^1-\frac{3i}{4} e^{z\mp3ix}(y^2-1)X^{2\pm}-\frac{3}{4}e^zyX^3,\\
Z^2&=\frac{1}{16y}(3qe^{\frac{z}{2}\pm\frac{3ix}{2}}(y\pm i) (5y^2\mp 12iy-5)\pm ie^{\frac{z}{2}\mp\frac{3ix}{2}} (y\mp i) (5y^2\pm 12i y-5))X^1\\
&\quad{}+\frac{i}{16y}(3qe^{\frac{z}{2}\mp\frac{3ix}{2}}(y \mp i) (5y^2\mp 2i y-5)\pm i e^{\frac{z}{2}\mp \frac{9ix}{2}}(y\pm i)  (5y^2\pm 2i y-5))X^{2\pm}\\
&\quad{}+\frac{1}{4}e^{\frac{z}{2}}(3qe^{\pm \frac{3ix}{2}}(y\mp i)\pm i e^{\mp \frac{3ix}{2}}(y \pm i))X^3,\\
Z^3&=\frac{3}{16y}e^{\frac{z}{2}\pm \frac{3ix}{2}}(5y^3\mp 7 i y^2+7 y\mp 5i)X^1+\frac{3i}{16y}e^{\frac{z}{2}\mp\frac{3ix}{2}}(5y^3\mp 7 i y^2-7y\pm 5i)X^{2\pm}\\
&\quad{}+\frac{3}{4}e^{\frac{z}{2}\pm \frac{3ix}{2}}(y\mp i)X^3.
\end{align*}

We obtain 
\begin{align}\label{s1b}
S^1&=\frac{1}{16y}e^{\frac{z}{2}}(3qe^{\pm\frac{3ix}{2}}(y\pm i)(y^2\mp 4i y-1)\pm i  e^{\mp\frac{3i x}{2}}(y\mp i)(y^2\pm 4i y-1))X^1\\\nonumber
&\quad{} + \frac{i}{16y}e^{\frac{z}{2}}(3q e^{\mp\frac{3i x}{2}}(y\mp i)(y^2\mp 2iy-1)\pm ie^{\mp\frac{9ix}{2}}(y\pm i)(y^2\pm 2iy-1))X^{2\pm},\\\label{s2b}
S^2&=\frac{3}{16y}e^{\frac{z}{2}\pm\frac{3ix}{2}}(y^3\mp 3 i y^2+3y\mp i)X^1+ \frac{3i}{16y}e^{\frac{z}{2}\mp \frac{3ix}{2}}(y^3\mp 3iy^2-3y\pm i)X^{2\pm},\\\label{s3b}
S^3&=\frac{2}{3y^2}e^{-z}(8y^2(y^2+1)p^2+4y(3y^2+1)p+5y^2+1)X^1\\\nonumber
&\quad{}+ i \frac{2}{3y^2}e^{-z\mp 3ix}(8 y^2(y^2-1)p^2+4y(3y^2-1)p+5y^2-1)X^{2\pm}\\
&\quad{}+\frac{4(2yp+1)^2}{3y}e^{-z}X^3.\nonumber
\end{align}
Observe that since $E^1$ and $E^3$ are imaginary, 
\[
X^3=\frac{1}{2y}(i(y^2-1) E^1+(y^2+1)E^2+10 i E^3)
\]
is real and also $\overline{i X^{2+}}=i X^{2-}$. It follows again that $S^{1-}$, $S^{2-}$, $S^{3-}$ are the complex conjugates of $S^{1+}$, $S^{2+}$, $S^{3+}$ respectively, depending on the sign of $\al$.

\begin{theorem}\label{g2b}
For the maximally symmetric distribution spanned by $X^1$, $X^{2\pm}$ when $\al=\pm 3i$ and annihilated by the 1-forms in (\ref{rolling1form}) with $(\si_1,\si_2,\si_3)$ given by (\ref{su2i}), the Lie algebra of split $\frak{g}_2$ is obtained from the pairwise bracket-generating set $\{S^1,S^2,S^3\}$ where $S^1$, $S^2$ and $S^3$ are given by (\ref{s1b}),  (\ref{s2b}),  (\ref{s3b}).
\end{theorem}

Using the transformation given by either
\begin{align*}
(x,y,z,p,q)\rightarrow \left(\frac{1}{\al'}\log(\frac{2}{9}e^z\al' y),\frac{1}{\al' \al y}, \log(\frac{\al' e^{\al x}y}{18}),-\frac{\al'}{18\al }y(\al qe^{\al x}-1),-\frac{1}{2y}(2yp+1)e^{-z}\right)
\end{align*}
or
\begin{align*}
(x,y,z,p,q)\rightarrow \left(\frac{1}{\al'}\log(-\frac{2}{9}e^z\al' y),\frac{9 \al y}{\al'}, \log(\frac{\al' e^{\al x}}{18y}),-\frac{\al'}{18\al y}(\al qe^{\al x}+1),-\frac{1}{2y}(2yp+1)e^{-z}\right),
\end{align*}
we can map the rolling distribution (\ref{rolling1form}) with the complex 1-forms given in (\ref{su2i}) from the parameter $\al^2=-\frac{1}{9}$ to $\al'^2=-9$. 

Substituting $y=i e^{-i \psi}$, $z=-i(\theta-\psi)$, $p=-\frac{1}{2}(\phi e^{-i \theta}-i) e^{i \psi}$
maps the 1-forms $\si_1$, $\si_2$, $\si_3$ in (\ref{su2i}) to 
\begin{align*}
\si_1&=\sin(\psi)\der \theta-\cos(\psi) i e^{-i \theta}\der \phi=\sin(\psi)\der \theta-\cos(\psi) \sin(\theta)\der \phi-i \cos(\psi)\cos(\theta) \der\phi,\\
\si_2&=\cos(\psi)\der \theta+\sin(\psi) i e^{-i \theta}\der \phi=\cos(\psi) \der \theta+\sin(\psi) \sin(\theta) \der \phi+i \sin(\psi)\cos(\theta) \der \phi,\\
\si_3&=-\der \psi-e^{-i \theta}\der \phi=-\der \psi-\cos(\theta) \der \phi+i \sin(\theta) \der \phi.
\end{align*}
The dual vector fields are given by 
\begin{align*}
E^1&=\sin(\psi)\partial_{\theta}-\cos(\psi)(\sin(\theta) \partial_{\phi}+i (\partial_\psi-\cos(\theta)\partial_{\phi})),\\
E^2&=\cos(\psi)\partial_{\theta}+\sin(\psi)(\sin(\theta) \partial_{\phi}+i (\partial_\psi-\cos(\theta)\partial_{\phi})),\\
E^3&=-\partial_{\psi}.
\end{align*}
Under this change of coordinates $(y,z,p)$ are no longer real variables. Taking $\al=\frac{i}{3}$ for instance, we get the bracket-generating set 
\begin{align*}
S^{1+}&=\frac{3}{4}e^{\frac{i \theta}{2}+\frac{i x}{6}}\sin\left(\frac{\psi}{2}\right)X^1-\frac{3}{2}e^{\frac{i \theta}{2}-\frac{i x}{6}}\cos\left(\frac{\psi}{2}\right)X^{2+}+\frac{\phi}{2}S^{2+}\\
S^{2+}&=\frac{3}{2}e^{-\frac{i \theta}{2}+\frac{i x}{6}}\cos\left(\frac{\psi}{2}\right)X^1+3e^{-\frac{i \theta}{2}-\frac{i x}{6}}\sin\left(\frac{\psi}{2}\right)X^{2+},\\
S^{3+}&=-4iqX^1+\frac{2}{3}i (q^2-9e^{-\frac{2i x}{3}})X^{2+}-\frac{2}{9}q^2e^{\frac{ix}{3}}(5E^3+3iE^1)
\end{align*}
of Theorem \ref{g2a}.
When $\al=3i$, if we denote by $V^1$, $V^2$ the vector fields
\begin{align*}
V^1&=\sin(\psi) X^1-\cos(\psi)e^{-3ix}X^{2+},\\
V^2&=\cos(\psi) X^1+\sin(\psi)e^{-3ix}X^{2+},
\end{align*}
then we obtain
\begin{align*}
S^{1+}&=-\frac{1}{8}e^{-\frac{i}{2}(3x+\theta)}\bigg(\cos\left(\frac{\psi}{2}\right)V^1-\frac{1}{2}\sin\left(\frac{\psi}{2}\right)V^2\bigg)+qS^{2+},\\
S^{2+}&=\frac{3}{2}e^{\frac{i}{2}(3x-\theta)}\bigg(\sin\left(\frac{\psi}{2}\right)V^1+\frac{1}{2}\cos\left(\frac{\psi}{2}\right)V^{2}\bigg),\\
S^{3+}&=-\frac{4}{3}ie^{i \theta}V^1-\frac{8}{3}i \phi V^2+\frac{4}{3}i e^{-i \theta}\phi^2(2V^1-\cos(\psi) E^1+\sin(\psi)E^2+5i E^3)
\end{align*}
of Theorem \ref{g2b}.
The action of the complex vector fields $(E^1, E^2, E^3)$ on the coordinate functions $(\theta, \phi, \psi)$ is quite different from that given by the real left-invariant vector fields on $SU(2)$.
If we take the complex conjugate of $\{\si_1, \si_2, \si_3\}$, i.e.\ take the complex parametrisation 
\begin{align*}
\bar\si_1&=-i(\der y+y\der z)+i(y^2-1)(\der p-p \der z)=-\si_1,\nonumber\\
\bar \si_2&=\der y+y\der z-(y^2+1)(\der p-p \der z)=\si_2, \\
\bar\si_3&=2iy(\der p-p \der z)-i \der z=-\si_3, \nonumber
\end{align*} 
and make the substitition  $y=-i e^{i \psi}$, $z=i(\theta-\psi)$, $p=-\frac{1}{2}(\phi e^{i \theta}+i) e^{-i \psi}$, then we obtain the 1-forms
\begin{align*}
\bar \si_1&=\sin(\psi)\der \theta+\cos(\psi) i e^{i \theta}\der \phi,\\
\bar \si_2&=\cos(\psi)\der \theta-\sin(\psi) i e^{i \theta}\der \phi,\\
\bar \si_3&=-\der \psi-e^{i \theta}\der \phi.
\end{align*}
To relate the complex parametrisation of $su(2)$ with the real rolling case, we see that the left-invariant 1-forms on $SU(2)$ can be obtained from the real part of $\si_1$, $\si_2$ and $\si_3$. The real 1-forms of the rolling distribution obtained in \cite{r17} can then be recovered by taking the real part $\frac{1}{2}(\om_1+\bar \om_1)$,  $\frac{1}{2}(\om_2+\bar \om_2)$,  $\frac{1}{2}(\om_3+\bar \om_3)$ and the corresponding real vector fields be computed from the complex distribution. It would be interesting to see if it is possible to write down the analogous bracket-generating set of vector fields for split $\frak{g}_2$ in the real case of the rolling distribution. It is also worthwhile to see if the vector fields $S^1$, $S^2$ and $S^3$ presented in Theorems (\ref{g2a}) and (\ref{g2b}) can be presented in a more symmetric fashion, following \cite{r17}.

\qr 

\end{document}